\documentclass[12pt]{article}
\usepackage{amssymb}
\usepackage{amsmath}
\usepackage[mathscr]{eucal}

\begin{document}

\def\cC{{\cal C}}
\def\cF{{\cal F}}
\def\cH{{\cal H}}
\def\cP{{\cal P}}
\def\bd{\mathop{\rm bd}}
\def\cl{\mathop{\rm cl}}
\def\rdiv{\mathop{\rm div}}
\def\const{\mathop{\rm const}}
\def\mn{{\medskip\noindent}}

\newtheorem{theorem}{Theorem}{\bf}{\it}
\newenvironment{form}{\begin{list}{}{\labelwidth 1.1cm
\labelsep 0.2cm \leftmargin 1.3cm}}{\end{list}}

\begin{center}
{\Large\bf An Arithmetic Proof of John's Ellipsoid Theorem}
\end{center}

\begin{center}
Peter M. Gruber and Franz E. Schuster
\end{center}

\begin{quote}
\footnotesize{ \vskip 1truecm\noindent {\bf Abstract.} Using an
idea of Voronoi in the geometric theory of positive definite
quadratic forms, we give a transparent proof of John's
characterization of the unique ellipsoid of maximum volume
contained in a convex body. The same idea applies to the \lq hard
part\rq\ of a generalization of John's theorem and shows the
difficulties of the corresponding \lq easy part\rq .

\medskip
\noindent {\bf MSC 2000.}\ 52A21, 52A27, 46B07.

\medskip\noindent
{\bf Key words.} John's theorem, approximation by ellipsoids,
Banach--Mazur distance.}
\end{quote}

\vskip 1truecm \centerline{\bf Introduction and Statement of
Results}

\vskip 0.5truecm\noindent The following well-known
characterizations of the unique ellipsoid of maximum volume in a
convex body in Euclidean $d$-space are due to John \cite{John48}
((i)$\Rightarrow$(ii)) and Pelczy\'nski \cite{Pelczynski90} and
Ball \cite{Ball92} ((ii)$\Rightarrow$(i)), respectively. For
references to other proofs, a generalization and to the numerous
applications see \cite{Ball01,GiannopoulosMilman01,JohnLind01}.

\begin{theorem}
Let $C\subset\mathbb E^d$ be compact, convex, symmetric in the
origin $o$, and with $B^d\subset C$. Then the following claims are
equivalent:
\partopsep -6pt
\parsep -6pt
\begin{itemize}
\itemsep-3pt
\item[{\rm (i)}] $B^d$ is the unique ellipsoid of maximum volume
in $C$.
\item[{\rm (ii)}] There are $u_k\in B^ d\cap \bd C$ and
$\lambda_k>0, k=1,\dots ,n,$ where $d\le n\le {1\over 2} d(d+1)$,
such that
$$I=\sum\limits_k \lambda_k\, u_k\otimes u_k.$$
\end{itemize}
\end{theorem}

\noindent Here, $B^d$ is the solid unit ball in $\mathbb E^d$, $I$
the $d\times d$ unit matrix, and for $u,v\in\mathbb E^d$ the
$d\times d$ matrix $u\, v^T$ is denoted by $u\otimes v$. $\bd$
stands for boundary.

\begin{theorem}
Let $C\subset\mathbb E^d$ be compact, convex, and with $B^d\subset
C$. Then the following claims are equivalent:

\partopsep -6pt
\parsep -6pt
\begin{itemize}
\itemsep-3pt
\item[{\rm (i)}] $B^d$ is the unique ellipsoid of maximum volume
in $C$.
\item[{\rm (ii)}] There are $u_k\in B^ d\cap \bd C$ and
$\lambda_k>0, k=1,\dots ,n,$ where $d+1\le n\le {1\over 2}
d(d+3)$, such that
$$I=\sum\limits_k \lambda_k\, u_k\otimes u_k,\, o=\sum\limits_k
\lambda_k\, u_k.$$
\end{itemize}
\end{theorem}

Our proof of Theorem 1 is based on the idea of Voronoi in the
geometric theory of positive definite quadratic forms, to
represent ellipsoids in $\mathbb E^d$ with center $o$ by points
in $\mathbb E^{{1\over 2}d(d+1)}$, see
\cite{EGH,Gruber05,RyshkovBaranovskii}. The problem on maximum
volume ellipsoids in $\mathbb E^d$ is then transformed into a
simple problem on normal cones in $\mathbb E^{{1\over 2}d(d+1)}$,
which can be solved easily by Carath\'eodory's theorem on convex
hulls. This idea has been applied before by the first author
\cite{ Gru88}. The proof of Theorem 2 is a simple extension. The
proof of the latter also gives Theorem 4 of Bastero and Romance
\cite{BasteroRomance}, where $B^d$ is replaced by a compact
connected set with positive measure.

In the context of John's theorem, it is natural to ask whether
ellipsoids can be replaced by more general convex or non-convex
sets. The following is a slight refinement of results of
Giannopoulos, Perissinaki and Tsolomitis \cite{GianPerTsol01} and
Bastero and Romance \cite{BasteroRomance} (Theorem 3). The result
of Giannopoulos et. al. was first observed by Milman in the case,
where both bodies are centrally symmetric, see \cite{Tomczak}.

\begin{theorem}
Let $C\subset\mathbb E^d$ be compact and convex, and $B\subset C$
compact with positive measure. Then {\rm (i)} implies {\rm (ii)},
where the claims {\rm (i)} and {\rm(ii)} are as follows:

\partopsep -6pt
\parsep -6pt
\begin{itemize}
\itemsep-3pt
\item[{\rm (i)}] $B$ has maximum measure amongst all its affine
images contained in $C$.
\item[{\rm (ii)}] There are $u_k\in B\cap \bd C, v_k\in N_C(u_k)$,
and $\lambda_k>0, k=1,\dots , n,$ where $d+1\le n\le d\,(d+1)$,
such that
$$I=\sum\limits_k \lambda_k \,u_k\otimes v_k,\,\, o=\sum\limits_k
\lambda_k\, v_k.$$  \end{itemize}
\end{theorem}
\medskip
\noindent Here $N_C(u), u\in\bd C$, is the normal cone of $C$ at
$u$. For this concept and other required notions and results of
convex geometry we refer to \cite{Schneider}.

Note that $B$ is not necessarily unique. A suitable modification
of Voronoi's idea applies in the present context and thus leads
to a proof of Theorem 3, paralleling our proofs of Theorems 1 and
2. Incidentally, the proof of Theorem 3 shows, why it is {\it not}
clear that property (ii) implies property (i), see the Final
Remarks.

\vskip 1truecm\centerline{\bf Proof of Theorem~1}

\vskip 0.5truecm\noindent For (real) $d\times d$--matrices
$A=(a_{ij}), B=(b_{ij})$ define $A\cdot B=\sum a_{ij} b_{ij}$.
The dot $\cdot$ denotes also the inner product in $\mathbb E^d$.
Easy arguments yield the following:

\begin{form}
\item[\hbox to 1.1cm{\rm (1)\hfill}]
Let $M$ be a $d\times d$ matrix and $u,v,w\in\mathbb E^d$. Then $M
u\cdot v=M\cdot u\otimes v$ and $(u\otimes v)w=(v\cdot w)u$.
\end{form}
Next, we specify two tools:

\begin{form}
\item[\hbox to 1.1cm{\rm (2)\hfill}]
Each $d\times d$ matrix $M$ with $\det M\not= 0$ can be
represented in the form $M=A R$, where $A$ is a symmetric,
positive definite and $R$ is an orthogonal $d\times d$ matrix.
\end{form}
(Put $A=(MM^T)^{1\over 2}, R=A^{-1}M$, see \cite{Gelfand67},
p.112.) Identify a symmetric $d\times d$ matrix $A=(a_{ij})$ with
the point $(a_{11},\dots , a_{1d}, a_{22},\dots , a_{2d}, \dots ,
a_{dd})^T \in \mathbb E^{{1\over 2}d(d+1)}$. The set of all
symmetric, positive definite $d\times d$ matrices then is
(represented by) an open convex cone $\mathscr P\subset\mathbb
E^{{1\over 2}d(d+1)}$ with apex at the origin. The set

\begin{form}
\item[\hbox to 1.1cm{\rm (3)\hfill}]
$\mathscr D=\{ A\in\mathscr P:\det A\ge 1\}$ is a closed, smooth,
strictly convex set in $\mathscr P$ with non-empty interior.
\end{form}
(Use the implicit function theorem and Minkowski's inequality for
symmetric, positive definite $d\times d$ matrices, see
\cite{RobertsVarberg}, p.205.)

(i)$\Rightarrow$(ii): By (2), any ellipsoid in $\mathbb E^d$ can
be represented in the form $A B^d$, where $A\in\mathscr P$. Thus
the family of all ellipsoids in $C$ is represented by the set
$$\mathscr E=\{ A\in\mathscr P:A u\cdot
v=A\cdot u\otimes v\le h_C(v) \text{ for } u,v\in S^{d-1}\},$$
see (1). Here, $h_C(\cdot )$ is the support function of $C$.
Clearly, $\mathscr E$ is the intersection of the closed halfspaces

\begin{form}
\item[\hbox to 1.1cm{\rm (4)\hfill}]
$\{ A\in\mathbb E^{{1\over 2}d(d+1)}:A\cdot u\otimes v\le
h_C(v)\}:u,v\in S^{d-1},$
\end{form}
with the set $\mathscr P$. Thus, in particular, $\mathscr E$ is
convex. By (i), $\mathscr E\setminus\{I\}\subset\{A\in\mathscr P:
\det A<1\}$. This, together with (3), shows that

\begin{form}
\item[\hbox to 1.1cm{\rm (5)\hfill}]
$\mathscr D$ and $\mathscr E$ are convex, $\mathscr D\cap\mathscr
E=\{I\}$, and $\mathscr D$ and $\mathscr E$ are separated by the
unique support hyperplane $\mathscr H$ of $\mathscr D$ at $I$ in
$\mathbb E^{{1\over 2}d(d+1)}$.
\end{form}

$\mathscr E$ is the intersection of the closed halfspaces in (4)
with the set $\mathscr P$, and these halfspaces vary continuously
as $u,v$ range over $S^{d-1}$. Thus the support cone $\mathscr K$
of $\mathscr E$ at $I$ can be represented as the intersection of
those halfspaces, which contain $I$ on their boundary
hyperplanes, i.e.\ for which $I\cdot u\otimes v=u\cdot v=h_C(v)$.
Since $u\cdot v\le 1$ and $h_C(v)\ge 1$ and equality holds in
both cases precisely when $u=v\in S^{d-1}\cap\bd C$ (note that
$B^d\subset C)$, we see that

\begin{form}
\item[\hbox to 1.1cm{\rm (6)\hfill}]
$\displaystyle \mathscr K=\bigcap\limits_{u\in B^d\cap \bd C}
\{A\in \mathbb E^{{1\over 2}d(d+1)}:A\cdot u\otimes u\le 1\} .$
\end{form}
The normal cone $\mathscr N$ of ($\mathscr E$ or) $\mathscr K$ at
$I$ is generated by the exterior normals of these halfspaces,

\begin{form}
\item[\hbox to 1.1cm{\rm (7)\hfill}]
$\mathscr N=\text{pos}\, \{u\otimes u:u\in B^d\cap \bd C\}.$
\end{form}

The cone $\mathscr K$ has apex $I$ and, by (5), is separated from
the convex set $\mathscr D$ by the hyperplane $\mathscr H$, where
$\mathscr H$ is the unique support hyperplane of $\mathscr D$ at
$I$. By considering the gradient of the function $A\to\det A$, we
see that $I$ is an interior normal vector of $\mathscr D$ at $I$
and thus a normal vector of $\mathscr H$ pointing away from
$\mathscr K$. Hence $I\in\mathscr N$. (7) and Carath\'eodory's
theorem for convex cones then yield the following: there are
$u_k\otimes u_k\in \mathscr N$, i.e. $u_k\in B^d\cap \bd C,$ and
$\lambda_k>0$ for $ k=1,\dots , n,$ where $n\le {1\over
2}d(d+1)$, such that
\begin{form}
\item[\hbox to 1.1cm{\rm (8)\hfill}]
$I=\sum\limits_k \lambda_k \,u_k\otimes u_k.$
\end{form}
For the proof that $n\ge d$, it is sufficient to show that
$\text{lin}\, \{ u_1,\dots , u_n\}=\mathbb E^d$. If this were not
true, we could choose $u\not= o, u\perp u_1,\dots , u_n$, and then
(1) yields the contradiction
$$0\not= u^2=Iu\cdot u=(\sum\limits_k \lambda_k\,(u_k\otimes
u_k)\,u)\cdot u=(\sum\limits_k \lambda_k \,(u_k\cdot u)\,u_k)\cdot
u=0.$$

(ii)$\Rightarrow$(i): Let $\mathscr E$ be as above. $\mathscr E$
is convex. $B^d\subset C$ implies that $I$ satisfies all defining
inequalities of $\mathscr E$, in particular those corresponding
to $u=v=u_k, k=1,\dots , n$. Since $h_C(u_k)=1$, these
inequalities are satisfied even with the equality sign. Thus
$I\in\bd \mathscr E$. Define $\mathscr K, \mathscr N$ and
$\mathscr H$ as before. (ii) implies that $I\in\mathscr N$. Hence
$\mathscr K$ is contained in the closed halfspace with boundary
hyperplane $\mathscr H$ through $I$ and exterior normal vector
$I$. Clearly, $\mathscr H$ separates $\mathscr K$ and $\mathscr D$
and thus, a fortiori, $\mathscr E (\subset \mathscr K)$ and
$\mathscr D$. Since $\mathscr D$ is strictly convex by (3), $
\mathscr D\cap \mathscr E = \{I\}$. Hence $B^d$ is the unique
ellipsoid of maximum volume in $C$.

\vskip 1truecm\centerline{\bf Outline of the Proof of Theorem~2}

\vskip 0.5truecm \noindent The proof of Theorem 2 is almost
identical with that of Theorem 1: an ellipsoid now has the form
$AB^d+a$ and is represented by $(A,a)\in \mathscr P\times \mathbb
E^d\subset\mathbb E^{\frac{1}{2}d(d+3)}$. $\mathscr E$ is the set
$$\{(A,a)\in\mathscr P\times \mathbb E^d:A\cdot u\otimes v+a\cdot
v\le h_C(v) \,\,\text{for}\,\,u,v\in S^{d-1}\}$$ and instead of
(5) we have
\begin{form}
\item[\hbox to 1.1cm{\rm \hfill}]
$\mathscr D\times \mathbb E^d$ and $\mathscr E$ are convex,
$(\mathscr D\times \mathbb E^d)\cap\mathscr E=\{(I,o)\}$ and
$\mathscr D\times \mathbb E^d$ and $\mathscr E$ are separated by
the hyperplane $\mathscr H\times \mathbb E^d$, where $\mathscr H$
is the unique support hyperplane of $\mathscr D$ at $I$ (in
$\mathbb E^{{1\over 2}d(d+1)}$).
\end{form}
$\mathscr K$ and $\mathscr N$ are the cones $$\displaystyle
\mathscr K=\bigcap\limits_{u\in B^d\cap \bd C}( A,a)\in \mathbb
E^{{1\over 2}d(d+3)}:A\cdot u\otimes u +a\cdot u\le 1\},$$
$$\mathscr
N=\text{pos}\, \{(u\otimes u,u):u\in B^d\cap \bd C\}.$$ As
before, $(I,o)\in\mathscr N$. Carath\'eodory's theorem for cones
in $E^{\frac{1}{2}d(d+3)}$ then shows the following: there are
$(u_k\otimes u_k, u_k)\in \mathscr N$ or, equivalently, $u_k\in
B^d\cap \bd C$ and $ \lambda_k>0,\,\, k=1,\dots , n,$ where $n\le
{1\over 2}d(d+3)$, such that instead of (8) we have
$$(I,o)=(\sum\limits_k \lambda_k \,u_k\otimes u_k,\sum_k\lambda_k
u_k).$$Since $o=\sum\lambda_k u_k$ and $\lambda_k>0$, the proof
that $n\ge d+1$ is the same as that for $n\ge d$ above. This
concludes the proof that (i)$\Rightarrow$(ii). The proof of
(ii)$\Rightarrow$(i) is almost the same as that of the
corresponding part of the proof of Theorem 1.

\vskip 1truecm\centerline{\bf Proof of Theorem~3}

\vskip 0.5truecm \noindent Identify a \,$d\times d$ matrix
$M=(m_{ij})$ with the point $(m_{11}, \dots , m_{1d},$ $m_{21},
\dots ,$ $m_{2d}, \dots ,$ $m_{dd})^T\in\mathbb E^{d^2}$. The set
$\mathscr P^\prime$ of all non-singular $d\times d$ matrices then
is (represented by) an open cone in $\mathbb E^{d^2}$ with apex
at the origin. The set

\begin{form}
\item[\hbox to 1.1cm{\rm \hfill}]
$\mathscr D^\prime = \{ M\in\mathscr P^\prime : |\det M|\ge 1\}$
is a closed body in $\mathscr P^\prime$, i.e.\ it is the closure
of its interior, with a smooth boundary surface.
\end{form}
The set of all affine images of $B$ in $C$ is represented by the
set
$$\mathscr E^\prime = \{ (M,a)\subset\mathscr P^\prime \times
\mathbb E^d: M u\cdot v+a\cdot v=M\cdot u\otimes v+a\cdot v\le
h_C(v) \text{ for } u\in B, v\in S^{d-1}\}.$$ This set is the
intersection of the closed halfspaces

\begin{form}
\item[\hbox to 1.1cm{\rm (9)\hfill}]
$\{ (M,a)\in\mathbb E^{d(d+1)}:M\cdot u\otimes v+a\cdot v\le
h_C(v)\}:u\in B, v\in S^{d-1},$
\end{form}
and thus of a convex set, with the set $\mathscr P^\prime \times
\mathbb E^d$. Choose a convex neighborhood $\mathscr U^\prime$ of
$(I,o)\, (\in\mathscr E^\prime \cap (\mathscr P^\prime \times
\mathbb E^d))$ which is so small that it is contained in the open
set $\mathscr P^\prime \times \mathbb E^d$. By (i),

\begin{form}
\item[\hbox to 1.1cm{\rm \hfill}]
the convex set $\mathscr E^\prime \cap \mathscr U^\prime$ and the
smooth body $\mathscr D^\prime \times \mathbb E^d$ only have
boundary points in common, one being $(I,o)$.
\end{form}
Hence $\mathscr E^\prime \cap \mathscr U^\prime$ and thus the
support cone $\mathscr K^\prime$ of $\mathscr E^\prime \cap
\mathscr U^\prime$ at $(I,o)$ is contained in the closed
halfspace whose boundary hyperplane is the tangent hyperplane of
the smooth body $\mathscr D^\prime \times \mathbb E^d$ at $(I,o)$
and with exterior normal pointing into $\mathscr D^\prime \times
\mathbb E^d$. This normal is $(I,o)$. The normal cone $\mathscr
N^\prime$ of $\mathscr K^\prime$ thus contains $(I,o)$.

The support cone $\mathscr K^\prime$ is the intersection of those
halfspaces in (9), which contain the apex $(I,o)$ on their
boundary hyperplanes. Thus $I\cdot u\otimes v+o\cdot v=h_C(v)$,
which is  equivalent to $u\in B\cap \bd C, v\in N_C(u)$. Hence,
these halfspaces have the form
$$\{ (M,a)\in \mathbb E^{d(d+1)}:M\cdot u\otimes v+a\cdot v\le
h_C(v)\}:u\in B\cap\bd C,\, v\in N_C(u),$$where $N_C(u)$ is the
normal cone of $C$ at the boundary point $u$. Thus, being the
normal cone of $\mathscr K^{\prime}$,
$$\mathscr N^\prime = \text{pos}\, \{( u\otimes v, v):u\in
B\cap\bd C,\, v\in N_C(u)\}.$$ Since $(I,o)\in \mathscr N^\prime$,
Carath\'eodory's theorem for convex cones in $\mathbb E^{d(d+1)}$
yields the following: there are $(u_k\otimes v_k, v_k)\in \mathscr
N^\prime$ or, equivalently, $u_k\in B\cap\bd C, v_k\in N_C(u_k)$,
and $ \lambda_k>0,\,\, k=1,\dots , n,$ where $n\le d(d+1)$, such
that
$$(I,o)=(\sum\limits_k \lambda_k \,u_k\otimes v_k,\sum\limits_k
\lambda_k \,v_k).$$ For the proof that $n\ge d+1$ we show by
contradiction that ${\rm lin}\{v_1,\dots , v_n\}=\mathbb E^d$ as
in the proof of Theorem~1.

\vskip 1truecm\centerline{\bf Final Remarks}

\vskip 0.5truecm \noindent In different versions of the proofs of
Theorems 1 and 2, which are closer to Voronoi's idea, ellipsoids
are represented in the form $x^T Ax\le 1$ and $(x-a)^T A(x-a)\le
1$, respectively.

If in Theorem 3 claim (ii) holds, then the support cone $\mathscr
K^\prime$ of $\mathscr E^\prime \cap \mathscr U^\prime$ at
$(I,o)$ is contained in the halfspace whose boundary is the
tangent hyperplane of $\mathscr D^\prime \times \mathbb E^d$ at
$(I,o)$ and with exterior normal pointing into $\mathscr D^\prime
\times \mathbb E^d$. Since $\mathscr D^\prime \times \mathbb E^d$
is not convex, this does {\it not} guarantee that $\mathscr
D^\prime \times \mathbb E^d$ and $\mathscr E^\prime$ do not
overlap , i.e.\ that (i) holds.

\medskip\noindent{\bf Acknowledgements.} The work of
the second author was supported by the Austrian Science Fund
(FWF). For his helpful remarks we are obliged to Professor Paul
Goodey.

\hfill\parbox[t]{6truecm}{ Forschungsgruppe \hfill\par Konvexe
und Diskrete Geometrie \hfill\par Technische Universit\"at
Wien\hfill\par Wiedner Hauptstra\ss e 8--10/1046\hfill\par A--1040
Vienna, Austria\hfill\par peter.gruber@tuwien.ac.at\hfill\par
franz.schuster@tuwien.ac.at\hfill}

\end{document}